\documentstyle{article}%%[12pt]
\def\Bbb{\bf}

\def\Z{{\Bbb Z}}

\def\be{\begin{equation}}
\def\ee{\end{equation}}
\def\bea{\begin{eqnarray*}}
\def\eea{\end{eqnarray*}}
\newtheorem{defn}{Definition}

\newtheorem{thm}{Theorem}

\newtheorem{cor}{Corollary}
\newenvironment{xpl}{\mbox{ }\\ \\{\bf Example}\mbox{ }}{
\hfill $\Box$\mbox{}\bigskip}

\newenvironment{proof}{\medskip {\bf Proof.}}{\hfill \rule{.5em}{1em} \\}
\begin{document}
\sloppy
\title{Polarized 4-Manifolds,\\ Extremal K\"ahler Metrics, and\\
Seiberg-Witten Theory}

\author{ Claude LeBrun\thanks{Supported
in part by  NSF grant DMS-9204093.}
\\
SUNY Stony
 Brook
  }

\date{May, 1995}
\maketitle

\begin{abstract}  Using  Seiberg-Witten theory, it is shown
 that any
K\"ahler metric
of constant negative scalar curvature on a compact 4-manifold
$M$ minimizes the $L^2$-norm
of scalar curvature among   Riemannian metrics
compatible with a  fixed decomposition $H^2(M)=H^+\oplus H^-$.
This implies, for example, that any such metric
on a minimal ruled surface must be locally symmetric.
 \end{abstract}
%\vfill
%\pagebreak

\section{Introduction}

In the late 1950's, Calabi first posed the problem of
representing each K\"ahler class on a compact
complex manifold by a K\"ahler metric of
constant scalar curvature. This eventually led him
\cite{cal}  to define extremal K\"ahler metrics,
which  minimize the
functional $\int s^2 ~ d\mu$ over
 a fixed K\"ahler class; here $s$ denotes
the scalar curvature, and $d\mu$ denotes the metric
volume measure. Any K\"ahler metric of constant scalar
curvature is  extremal in this sense,
  but Calabi showed  by example that the converse
is generally false.

In real dimension 4,   new insights into
this problem can be gained by temporarily
 venturing outside the K\"ahlerian arena, and instead
working in a broader   Riemannian context.
 Instead of fixing a
K\"ahler class, we will fix a closely-related direct sum
decomposition  $H^2(M,{\Bbb R})=H^+\oplus H^-$.
Seiberg-Witten theory will then allow us to
see that any K\"ahler metric of constant
negative scalar curvature is an absolute
minimum of $\int s^2 d\mu$ among metrics compatible
with such a decomposition. As an application,
we will then classify  K\"ahler metrics
of constant negative scalar curvature on
 minimal ruled surfaces.

\section{Polarizations}
 \begin{defn}
Let $M$ be a smooth compact oriented 4-manifold.
A {\em polarization} of $M$ is a
maximal
linear subspace
$H^+\subset H^2(M, {\Bbb  R})$  for which
the restriction of the intersection form
 is positive-definite.
\end{defn}
Because the intersection
form is non-degenerate, every polarization
determines an orthogonal complement
$H^-$ with respect to the intersection form, and the intersection form
is negative-definite on this orthogonal complement; this puts
 polarizations of $M$
and   of the   reverse-oriented
manifold $\overline{M}$
  in   natural
one-to-one correspondence.
The dimensions $b_\pm :=\dim H^\pm$   are
important  homeomorphism invariants
   of $M$, and their difference
 $\tau = b_+-b_-$  is called the signature.
Given a polarization, we will
routinely invoke the decomposition
$$H^2(M)= H^+\oplus H^-$$
to uniquely express  elements $\alpha \in H^2$
as    $\alpha = \alpha^++\alpha^-$,
where $\alpha^\pm\in H^\pm$.

While the imposition of  a polarization  may seem
 frivolous, such choices
  arise quite naturally in Riemannian geometry.
Indeed, if $g$ is a smooth Riemannian
metric on $M$, then the space
$$H^+(g) := \{ [\varphi ]\in H^2(M) ~|~ \varphi \in C^\infty (\wedge^2),
{}~ d\varphi = 0 ,
{}~ \varphi = \star_g \varphi \}
$$
of self-dual $g$-harmonic 2-forms
is a polarization on $M$. We will say that
a Riemannian metric $g$ is {\em adapted}
to the polarization $H^+$ if $H^+(g)=H^+$.
A polarization will be called  a {\em metric
polarization} if there is at least one metric
adapted to it.

\begin{xpl}
Let $(M,J,g)$ be a compact K\"ahler manifold of complex
dimension 2, and let $\omega$ denote the associated K\"ahler form.
Let $\Re e H^{2,0}\subset H^2(M,{\Bbb R})$ denote the   de Rham
classes which are represented by real parts of holomorphic
2-forms. Then
$H^+(g)={\Bbb R}[\omega ] \oplus \Re e H^{2,0} .$
\end{xpl}

Algebraic geometers  sometimes use the
term ``polarization'' to denote a choice of
K\"ahler class $[\omega ]$ on a compact complex
manifold $(M,J)$. In light of the above example, our
  terminology may thus be
justified by the fact that the polarization $H^+(g)$
of   K\"ahler metric determines
the K\"ahler class  $[\omega]$
if the complex structure $J$ and
total volume $[\omega]^2/2$   are specified.

 Because the Hodge $\star$-operator is conformally invariant
on middle-dimensional forms,
the present notion of  polarization is
conformally invariant; that is, $H^+(g)=H^+(fg)$
for any smooth positive function $f$. Thus all our
conclusions about metrics adapted to a fixed polarization
will also imply results about global conformal invariants.

\section{Seiberg-Witten Theory}

Let $M$   be a smooth connected
compact oriented 4-manifold, and assume that
$M$ admits an orientation-compatible
almost-complex structure $J:TM\to TM$, $J^2=-1$.
Such an almost-complex structure
determines a spin$^c$ structure
on $M$, meaning a cohomology class
$c\in H^2(F, \Z)$ on the
oriented frame bundle $F\to M$
whose restriction of $c$ to a typical
fiber $F_x\cong SL(4,{\Bbb R})\times {\Bbb R}$
is the non-zero element of $H^2(F_x, \Z)\cong \Z_2$.
Indeed, if $g$ is any $J$-invariant Riemannian
metric on $M$, let $F_{SO} \subset F$ be the bundle of
oriented $g$-orthogonal frames, and let $F_{U} \subset F_{SO}$
denote the bundle of unitary frames with respect to $g$
and $J$. Then the Poincar\'e dual   of the submanifold
$F_{U} \subset F_{SO}$,
thought of as an element of $H^2(F, \Z )=H^2(F_{SO} , \Z)$,
is fiber-wise non-zero   and is independent of $g$; this is the
  promised   spin$^c$ structure $c$.
 If a spin$^c$ structure   arises in this
way, we will say that it is of {\em almost-complex type},
and we will say that the almost-complex structure $J$ and
the spin$^c$ structure $c$ are {\em compatible}.
One may choose to think of  a spin$^c$ structure of
almost-complex type
as an equivalence class of almost-complex structures $J$;
 two such structures are then equivalent iff their graphs are
homologous as submanifolds of the bundle $F/GL(2, {\Bbb C})$
of orientation-compatible
almost-complex structures.

Any spin$^c$ structure on $M$ determines, up to isomorphism,
 a complex line bundle
$L\to M$ such that $c_1(L)\equiv w_2(TM) \bmod 2 , $
by setting $c_1(L)=2c\in H^2(M,\Z)\subset H^2(F,\Z )$;
and conversely any such a choice of $c_1(L)\in H^2(M,\Z)$
determines a spin$^c$ structure up to
  2-torsion in $H^2(M,\Z)$.
If we choose a Riemannian metric $g$ on $M$,
a spin$^c$-structure determines rank-2 Hermitian vector bundles
$V_\pm\to M$ with $\wedge^2V_\pm =L$ and $T^*M\otimes {\Bbb  C}
\cong\mbox{Hom}(V_+,V_-)$;  and on any contractible open
set in $M$ we have   canonical  (sign-ambiguous) isomorphisms
$$V_\pm = {\Bbb S}_\pm \otimes L^{1/2} , $$
where ${\Bbb S}_{\pm}$ are the left- and right-handed
spinor bundles of $g$, and $L^{1/2}$ is a complex
line bundle whose square is $L$. Each unitary
connection $A$ on $L$
therefore   induces a unitary  connection
$\nabla_A : C^{\infty}(V_+)\to C^{\infty}(T^*M\oplus V_+)$
on $V_+$, and following this with the isomorphism
$T^*M\otimes {\Bbb  C}
\cong\mbox{Hom}(V_+,V_-)$ gives us a {\em Dirac operator}
$D_A: C^{\infty}(V_+)\to C^{\infty}(V_-)$.

This can all be made much more concrete for spin$^c$ structures of
almost-complex type.
 Given a Riemannian metric $g$ on $M$, we can represent such a
spin$^c$ structure
by an almost-complex structure $J:TM\to TM$, $J^2=-1$ such
that $J^*g=g$.  The tangent bundle
$TM$ of $M$ is thereby given the structure of a rank-2 complex
vector bundle
$T^{1,0}$ by defining scalar multiplication by $i$ to be  $J$. Setting
$\wedge^{0,p}:=\wedge^p\overline{T^{1,0}}^*\cong \wedge^pT^{1,0}$,
the    bundles
$V_{\pm}$ of twisted spinors are given  by
\begin{eqnarray} V_+&=& \wedge^{0,0}\oplus \wedge^{0,2}\label{spl}\\
V_-&=&\wedge^{0,1}, \label{spr}\end{eqnarray}
and their Hermitian structures are the obvious ones induced by $g$.
In particular, $L$ is the anti-canonical bundle of $(M,J)$, and
we therefore have
$$c_1(L)^2=(2\chi + 3\tau )(M).$$
Spin$^c$ structures of almost-complex type
are characterized by this  last property.

If $(M,g,J)$ is a K\"ahler manifold,
so that
$J$ is parallel with respect to the metric connection,
and if $A$ is the so-called Chern connection on the
anti-canonical bundle $L$, then the connection
$\nabla_A$ on $V_+$ has a parallel sections
corresponding to the constant sections of
$\wedge^{0,0}\subset V_+$. Conversely,
a metric is K\"ahler for
$c$-compatible complex structure $J$ provided there
is a choice of $A$ for which $V_+$ has a parallel
section; indeed, this   implies that the
holonomy of ${\Bbb S}_+$ is contained in $U(1)\subset SU(2)$,
and the Riemannian holonomy is therefore contained
in $(U(1)\times SU(2))/\Z_2 =U(2)$.
For $g$ a K\"ahler metric and $A$ the Chern connection,
the Dirac operator can correspondingly be expressed as
$D_A ={\sqrt{2}}(\overline{\partial}\oplus \overline{\partial}^*)$,
where  $\overline{\partial}: C^{\infty}(\wedge^{0,0})\to
C^{\infty}(\wedge^{0,1})$ is  the $J$-antilinear part of the exterior
differential $d$, acting on complex-valued functions, and where
$\overline{\partial}^*: C^{\infty}(\wedge^{0,2})\to
C^{\infty}(\wedge^{0,1})$ is the formal  adjoint of
the  map induced by the exterior
differential $d$ acting on 1-forms;
more generally, $D_A$ will differ from
${\sqrt{2}}(\overline{\partial}\oplus \overline{\partial}^*)$
by only $0^{th}$ order terms.

  Let us now fix a spin$^c$ structure  $c$
 of almost-complex type on
on $M$.  For each Riemannian metric $g$,
the Seiberg-Witten equations \cite{witten}
\begin{eqnarray} D_{A}\Phi &=&0\label{drc}\\
 F_{A}^+&=&i \sigma(\Phi).\label{sd}\end{eqnarray}
are then
equations
for an unknown smooth connection $A$ on $L$
and an unknown
smooth section $\Phi$ of $V_+$.
Here the purely imaginary 2-form $F_{A}^+$  is the self-dual part of
the curvature of $A$, and
the natural real-quadratic map $\sigma: V_+\to \wedge^2_+$
satisfies $|\sigma (\Phi)|^2=|\Phi|^4/8$.
  For our purposes, it is crucial  that
equations (\ref{drc}) and (\ref{sd})
imply the Weitzenb\"ock formula
\be\label{wb}
 \nabla_A^*\nabla_A \Phi + \frac{s+|\Phi|^2}{4}\Phi  =0.
\ee
Given a  solution $(A , \Phi)$ of (\ref{drc}) and (\ref{sd})
and a mooth map
 $f: M\to S^1\subset {\Bbb C}$, the pair  $(\hat{A},\hat{\Phi})=
(A - 2f^{-1}d f,  f\Phi)$ is also a solution;
   solutions which are
related in this way are called {\em gauge equivalent}. A solution
is called {\em reducible} if $\Phi\equiv 0$; otherwise, it is
called {\em irreducible}.

Now, in addition to such a  spin$^c$ structure $c$,
let us fix a metric polarization $H^+$
on $M$. Assume that  $c_1^+:=[c_1(L)]^+\in H^+$ is non-zero,
which guarantees that
every solution of the Seiberg-Witten equations is
irreducible whenever  $g$ is  an $H^+$-adapted metric.
For each such metric, one can then define the {\em Seiberg-Witten
invariant} $n_c(M,g)\in \Z$ by counting the
 solutions modulo gauge equivalence with appropriate multiplicities
for a generic small perturbation of (\ref{drc}) and (\ref{sd});
in particular, $n_c(M,g)\neq 0$ implies there is a
solution of (\ref{drc}) and (\ref{sd}). We now define
$n_c(M,H^+)=n_c(M,g)$ for any $H^+$-adapted metric $g$.
This  is   metric-independent because
the moduli spaces corresponding to different $H^+$-adapted metrics
are  cobordant as oriented 0-manifolds.
Indeed, when $b_+>1$, $n_c(M,H^+)$ is even independent of the
polarization.  By contrast, when $b_+=1$, the invariant  generally
jumps \cite{FM,KM} as $H^+$ passes though polarizations for which
 $c_1^+=0$.

The following result \cite{FM,leb2,taubes,witten} shows that
the invariant is non-trivial for many interesting polarized manifolds:

\begin{thm}\label{non}
Let $(M,J,g)$ be a compact K\"ahler  surface
for which the K\"ahler class
 $[\omega ]$ satisfies $c_1\cdot [\omega] < 0$.
Let $H^+=H^+(g)$ be the metric polarization, and let
$c$ be the canonical spin$^c$ structure of $(M,J)$.
Then    $n_c(M,H^+)=1$.
 \end{thm}

 The assumption that $c_1\cdot [\omega] < 0$ amounts to the
requirement that the scalar curvature $s$  of $(M,g)$ be
negative ``on average,'' by virtue of the
Gauss-Bonnet-type formula
 $\int_M s ~d\mu = 4\pi c_1\cdot [\omega ]$.
In fact, the proof of Theorem \ref{non} becoms  particularly
simple \cite{leb}  if the scalar curvature
 is assumed to be a
negative constant.

The following scalar-curvature inequality is the
crux of the the present note:

\begin{thm}\label{est}
Let $(M,H^+)$ be a  polarized smooth compact
oriented  4-manifold, and suppose that
there is
  a spin$^c$ structure $c$ of almost-complex type
 on $M$ for which the
 Seiberg-Witten invariant
is non-zero; let $c_1(L)\in H^2(M, {\Bbb R})$ denote
the anti-canonical class of this structure, and
  let $c_1^+\neq 0$ be its orthogonal projection to
$H^+$ with respect to the intersection form.  Then  every
$H^+$-adapted Riemannian metric $g$ satisfies
$$ \int_M s^2~d\mu \geq 32\pi^2(c_1^+)^2 , $$
with equality iff $g$ is  K\"ahler with respect to a
$c$-compatible
complex structure and
 has constant negative
scalar curvature.
 \end{thm}
\begin{proof}
For any given metric $g$  adapted to $H^+$, there must exist an
irreducible solution
of  (\ref{drc}) and  (\ref{sd}), since otherwise we would have
$n_c(M, H^+)=0$. Now the Weitzenb\"ock formula \ref{wb} tells
us that
$$0= \int_M \left( 4|\nabla \Phi|^2 +
s |\Phi|^2 + |\Phi|^4 \right)~d\mu ,$$
so that $\int  (-s)|\Phi|^2~d\mu  \geq \int  |\Phi|^4 ~d\mu$, with
equality iff $\Phi$ is parallel.
The Schwartz inequality therefore tells us that
$$\int_M  s^2~d\mu  \geq
\frac{\left(\int_M  (-s)|\Phi|^2~d\mu \right)^2}{\int_M
  |\Phi|^4 ~d\mu}
\geq  \int_M  |\Phi|^4 ~d\mu
{}~,$$
with equality  iff $\nabla \Phi =0$ and $s$ is constant.
On the other hand,
$|F_A^+|^2=|\sigma (\Phi )|^2= |\Phi|^4/8$, so this may be rewritten as
$$\int_M  s^2~d\mu
\geq  8\int_M |F_A^+|^2    ~d\mu  .$$
But now  letting $\varphi$ denote the harmonic representative of
the de Rham class
$[F_A]=2\pi c_1$, we have
\bea
\int_M|F_A^+|^2d \mu &=& \frac{1}{2}\int_M (| F_A^+|^2-|F_A^-|^2)d \mu
+ \frac{1}{2}\int_M (|F_A^+|^2+|F_A^-|^2)d \mu
\\&=&2\pi^2c_1(L)^2+\frac{1}{2}\int_M|F_A|^2d \mu
\\&\geq &2\pi^2c_1(L)^2+\frac{1}{2}\int_M|\varphi|^2d \mu
\\&=& \frac{1}{2}\int_M (|\varphi^+|^2-|\varphi^-|^2)d \mu
+ \frac{1}{2}\int_M (|\varphi^+|^2+|\varphi^-|^2)d \mu
\\&=&  \int_M|\varphi^+|^2d \mu \\&=&
4\pi^2 (c_1^+)^2
 \eea
because a harmonic form  minimizes   the $L^2$ norm among
closed forms in its deRham class.
Hence
$$ \int_M s^2~d\mu \geq 32\pi^2(c_1^+)^2 , $$
as claimed. Moreover, equality is achieved only if  $s$ is constant and
 $\nabla \Phi =0$, which implies that   $g$ is K\"ahler with
respect to a $c$-compatible
complex structure.

 Conversely, any K\"ahler  metric with constant scalar curvature
will saturate this bound, since  $\varphi^+=s\omega /4$
and $d\mu=|\omega^2/2|$
 for such a metric.
\end{proof}

Notice that the above inequality will hold, more generally,
for any metric and  spin$^c$ structure for which there is an
irreducible solution of the Seiberg-Witten
equations. Thus, while it is  also  possible  in principal
to  define Seiberg-Witten
invariants for spin$^c$ structures
 which are not of almost-complex
type,  these  can  often be shown to vanish by
producing metrics for which $\int s^2 d\mu$ is
sufficiently small.

The above result has a curious ramification for
 conformal geometry. Let $g$ be a K\"ahler
metric of constant negative scalar curvature on a compact
complex surface $M$. Let $[g]$ be the
conformal class of $g$, and let $[g]_-\subset [g]$ be the
open subset  of metrics of negative scalar curvature.
Then $g$ simultaneously {\em minimizes} $\int s^2 d\mu$ and
{\em  maximizes} $(\int s d\mu)^2/\int d\mu$, considered
as functionals on $[g]_-$!

\section{Ruled Surfaces}

While the best-understood obstructions to the existence of
K\"ahler metrics of constant scalar curvature  entail
the existence of non-trivial holomorphic vector fields,
a more subtle obstruction, related to Mumford
stability, was discovered by Burns and
de Bartolomeis \cite{bb}. While their result deals only with
$s\equiv 0$ metrics on minimal ruled surfaces,
its formulation is so elegant as to make it desirable to
put this isolated result in a
 wider context.
 We shall now take a small step in
this direction by showing that analogous
conclusions apply to
  K\"ahler metrics of
constant {\em negative} scalar curvature on
minimal ruled surfaces.

\begin{thm} \label{lsym}
Let $M$ be the total space of
an oriented 2-sphere bundle $M\to\Sigma$ over a
compact oriented surface.
For some complex structure $J$,  suppose that  there is a
K\"ahler metric
$g$   of constant negative scalar curvature
on $M$. Then   the universal
cover of $(M,g)$ is isometric to the product
$S^2\times {\cal H}^2$, where
the 2-sphere
and  hyperbolic 2-space are  endowed with   appropriate
constant multiples of their standard metrics.
\end{thm}
\begin{proof}
First observe that the structure group $\mbox{Diff}^+(S^2)$
of
$M\to \Sigma$ can be reduced to $SO(3)$, since the induced
 induced metric on each fiber
$S^2$ can be  conformally rescaled  to yield
a metric of curvature $+1$, and
the freedom in doing so is paramaterized by the contractible
space $SL(2,{\Bbb C})/SU(2)$.
As a consequence, $M$ admits a fiberwise antipodal map
$\psi : M\to M$. Since this is orientation-reversing,
we may define a spin$^c$ structure $\bar{c}$ on the
reverse-oriented manifold $\overline{M}$ by setting
$\bar{c}=\psi^*c$, where $c$ is the spin$^c$ structure
on $M$   induced by $J$.

Because $\psi$ reverses
orientation, $\psi^*:H^2(M)\to H^2(M) $
reverses the sign of the intersection form:
 $$\psi^*(\alpha)\cdot \psi^* (\beta) = -\alpha \cdot \beta
{}~~~~~ \forall \alpha, \beta \in H^2(M, {\Bbb R}).$$
Since $(\psi^*)^2=(\psi^2)^*=1$,
it follows that $\alpha \cdot \psi^* (\alpha )=0$.
But $H^2(M,{\Bbb R})$ is 2-dimensional, and
 $\psi^*$ therefore takes
 any 1-dimensional subspace to its  intersection-form-orthogonal
subspace. This shows that $\psi^*(H^+(g))= H^-(g)$.

Since we therefore know that
$(M,H^+(g),c)$ and $(\overline{M}, H^-(g),\bar{c})$
are isomorphic  as
polarized oriented 4-manifolds with spin$^c$ structure,
they have the same Seiberg-Witten invariant, and
  $(c_1^+)^2(M,H^+(g),c)=(c_1^+)^2(\overline{M},H^-(g),\bar{c})$.
But $g$ is a K\"ahler metric with constant
negative scalar curvature, implying that $n_c(M,H^+(g))=1$ and
$\int_Ms^2 ~ d\mu = 32\pi^2 (c_1^+)^2(M,H^+(g),c)$.  Hence
 $n_c(\overline{M}, H^-(g))=1$ and
$\int_{M}s^2 ~ d\mu = 32\pi^2
(c_1^+)^2(\overline{M},H^-(g),\bar{c})$,
too. By Theorem \ref{est}, the latter implies that
$g$ is K\"ahler with respect to a complex structure
$\tilde{J}$ compatible with $\bar{c}$, and hence
compatible with the orientation of $\overline{M}$.

Now the $g$-preserving linear maps $J$ and $\tilde{J}$ lie
in opposite factors of $SO(4)=(SU(2)\times SU(2))/{\Bbb Z}_2$,
and so commute. The endomorphism $J\tilde{J}=\tilde{J}J$
is therefore diagonalizable, with eigenvalues
$\pm 1$, and we have an eigenspace decomposition
$TM=L_1\oplus L_2$ of the tangent bundle
into rank-2 real vector bundles.
Since both complex structures are invariant under
parallel transport with respect to $g$, this
decomposition is parallel, and the universal cover
of $(M,g)$ is therefore the Riemannian product
$(X_1,g_1)\times (X_2, g_2)$
of a pair of complete simply connected surfaces.
Since the scalar curvature $s$ of $g$ is constant,
and since $s$ is the sum $s_1+s_2$ of the scalar  curvatures
of  $(X_1,g_1)$ and $(X_2, g_2)$, it follows that $s_1$ and $s_2$ must
 both be constant.
Now $\pi_3(\Sigma )$ is finite, so the
exact homotopy sequence
$$\cdots \to \pi_3 (\Sigma ) \to \pi_2 (S^2 ) \to \pi_2 (M)
\to \pi_2 (\Sigma ) \to \cdots $$
predicts that $\pi_2(M)\neq 0$; passing to the
universal cover, we thus have
$\pi_2( X_1)\times \pi_2(X_2)=
\pi_2( X_1\times X_2)=\pi_2(M)\neq 0$, and   at least one of the
simply-connected surfaces $X_j$ must therefore
 be a 2-sphere. On the other hand,
 $s=s_1+s_2 < 0$, so the other
factor   must   be hyperbolic.
 \end{proof}

With this result in hand, we can now
 solve the existence and uniqueness problems for
 K\"ahler metrics of constant negative
scalar curvature on minimal ruled surfaces.

\begin{thm}\label{poly}
 Let $E \to \Sigma$ be a rank-2 holomorphic vector bundle
 over a
compact complex curve, and let
  $(M,J)={\Bbb P}(E)$
 be  the total space
of the associated ${\Bbb CP}_1$-bundle.
 Let $[\omega ]$ be
a K\"ahler class on $M$ with $c_1\cdot [\omega ] < 0$.
Then $[\omega ]$ contains a  K\"ahler metric of constant
scalar curvature
iff $E$ is a polystable
 vector bundle. Moreover, when such a metric exists,
 it is unique modulo biholomorphisms of $(M,J)$.
\end{thm}
\begin{proof}
Let us begin by reminding the reader
that  $\int s ~ d\mu = 4\pi c_1\cdot [\omega ]$
 for any K\"ahler metric in $[\omega ]$,
so the   $c_1\cdot [\omega ] < 0$ hypothesis
exactly limits our discussion to K\"ahler metrics of
constant {\em negative} scalar curvature.
Note that K\"ahler classes with $c_1\cdot [\omega ] < 0$ will exist on
the ruled surface $M\to \Sigma$
 iff $\Sigma$ has genus  $\geq 2$.

Now recall that a vector bundle $E$ is
said to be {\em polystable} (or sometimes {\em quasi-stable})
if it is a semi-stable bundle of the form $E=\bigoplus_{j=1}^n E_j$,
where the   $E_j$ are
 stable vector bundles. A landmark
result of Narasimhan and
Seshadri \cite{NS} asserts that, for bundles over a
Riemann surface, polystabilty is equivalent to
the existence of a flat projective unitary
connection on ${\Bbb P}(E)$. If $E$ is polystable, we
thus have $(M,J)=\Sigma \times_{\rho} {\Bbb CP}_1$ for a
representation
  $\rho : \pi_1 (\Sigma )\to PSU (2)=SO(3)$ of the fundamental
group which is unique up to conjugation in $SO(3)$.
When this happens, local products of
constant curvature K\"aler metrics on $\Sigma$ and
   ${\Bbb CP}_1$ provide us with a
constant scalar curvature K\"ahler metric in each
K\"ahler class. Our task is  therefore
to show that any K\"ahler metric on
${\Bbb P}(E)$ with $s =\mbox{const} < 0$ is
necessarily of this form.

To this end, assume that $g$ is a K\"ahler metric
of constant negative scalar curvature on $(M,J)$, and recall
 that Theorem \ref{lsym}
tells us that the universal cover of $(M,g)$
must be a product $S^2\times {\cal H}^2$
of spaces of constant   curvature. Once the factors are
correctly oriented, moreover, the product complex structure
will necessarily agree with the lift of $J$ because
  the holonomy of $S^2\times {\cal H}^2$ is $U(1)\times U(1)$.
Now $\pi_1(M)=\pi_1(\Sigma )$   acts on
$S^2\times {\cal H}^2$ by holomorphic isometries,
and thus sends any holomorphic curve $S^2\times \{ \mbox{pt}\}$
to another curve of this form --- after all, these are the {\em only}
compact complex curves in ${\Bbb CP}_1\times {\cal H}^2$! The
induced action on   ${\cal H}^2$ is, moreover,  free and proper,
since   $S^2$ is compact and every rotation of
$S^2$ has a fixed point. Thus
${\Bbb P}(E)=(S^2\times {\cal H}^2)/\pi_1(\Sigma)$ is
 biholomorphic to $\tilde{\Sigma}
\times_{\rho} {\Bbb CP}_1$ for  some compact
Riemann surface $\tilde{\Sigma}$  and  some representation
  $\rho : \pi_1 (\tilde{\Sigma} )\to PSU (2)=SO(3)$.
 But since the fibers of $\tilde{\Sigma}
\times_{\rho} {\Bbb CP}_1\to \tilde{\Sigma}$
 are
the only rational curves in $(M,J)$, there is a
biholomorphism  $\Sigma\to \tilde{\Sigma}$ such that the
diagram
$$\begin{array}{ccc}{\Bbb P}(E)
&\to & \tilde{\Sigma}\times_{\rho} {\Bbb CP}_1\\
\downarrow && \downarrow\\
\Sigma&\to & \tilde{\Sigma}
\end{array}$$
commutes. This gives $E$  a flat unitary projective connection,
and so shows that $E$ is polystable.
\end{proof}

The same reasoning can easily be applied to other
complex surfaces with orientation-reversing diffeomorphisms.
For example, on the product $\Sigma_1\times \Sigma_2$ of
Riemann surfaces of positive genus, any product  of
constant curvature metrics is the unique
constant-scalar-curvature metric in its K\"ahler
class. Closely related results have been proved by
Leung  \cite{leung}.

\begin{cor} Let $E \to \Sigma$ be a rank-2 holomorphic vector bundle
 over a
compact complex curve, and let
  $(M,J)={\Bbb P}(E)$.
 Let $[\omega ]$ be
a K\"ahler class on $M$ with $c_1\cdot [\omega ] < 0$.
If $[\omega ]$ contains an extremal  K\"ahler metric,
then   $E$ is either stable, or else  is the
direct sum $L_1\oplus L_2$ of a pair of holomorphic line bundles.
\end{cor}
\begin{proof}
By Theorem \ref{poly},
we may assume the scalar curvature $s$ of our extremal
metric $g$ is non-constant. Then  the isometry group of $g$
 is non-trivial, since $J\mbox{grad}_gs$ is
a Killing field. Hence there is a Killing field on
$M$ which generates a  non-trivial $U(1)$-action
by biholomorphisms. Since $\Sigma$ has
genus $\geq 2$, and so admits no non-trivial holomorphic
vector field, this action preserves the
fibers of $M\to \Sigma$, and has 2 distinct
fixed points in each fiber. The corresponding
linear subspaces of $E$ then give the desired
direct sum decomposition.
\end{proof}

\bigskip

{\noindent \bf Acknowledgement.} The author would like
 to thank Christina T{\o}nnesen for
some very useful conversations.

\end{document}